\documentclass[a4paper,11pt]{article}
\usepackage{mathtools}
\mathtoolsset{showonlyrefs}
\usepackage{amsmath}
\usepackage{amssymb}
\usepackage{amsthm}
\usepackage[utf8]{inputenc}
\usepackage[ruled]{algorithm2e}
\usepackage{etex}

\SetKw{KwInit}{Initialization:}
\SetKw{KwBreak}{break}
\SetKw{KwInput}{Input:}
\SetKw{KwOutput}{Output:}
\SetKw{KwIter}{Iteration:}

\usepackage{graphicx}
\usepackage{pgfplots}
\usepackage{url}

\makeatletter
\g@addto@macro\@floatboxreset\centering
\makeatother

\newcommand{\N}{\ensuremath{\mathbb{N}}}

\newcommand{\R}{\ensuremath{\mathbb{R}}}

\newcommand{\tT}{{\scriptscriptstyle \operatorname{T}}}
\DeclareMathOperator*{\argmin}{arg\,min}
\DeclareMathOperator{\prox}{prox}

\DeclareMathOperator*{\vect}{vec}

\usepackage{tikz,xargs}
\usetikzlibrary{spy}
\usetikzlibrary{positioning}

\title{
Removal of Curtaining Effects by a Variational Model with Directional Forward Differences
}

\newtheorem{remark}{Remark}
\newtheorem{proposition}{Proposition}

\author{Jan Henrik Fitschen\footnote{Department of Mathematics,
University of Kaiserslautern, Germany,
{fitschen@mathematik.uni-kl.de}}
\and Jianwei Ma\footnote{Department of Mathematics, Harbin Institute of 
Technology, Harbin, China, jma@hit.edu.cn}
\and Sebastian Schuff\footnote{Department of Mechanical and Process Engineering, 
University of Kaiserslautern, Germany,
{schuff@mv.uni-kl.de}} 
}

\begin{document}

\maketitle

\begin{abstract}
\noindent
Focused ion beam (FIB) tomography provides high resolution volumetric images on a micro scale.
However, due to the physical acquisition process the resulting images are often corrupted by a so-called curtaining or waterfall effect.
In this paper, a new convex variational model for removing such effects is 
proposed.
More precisely, an infimal convolution model is applied to split the corrupted 3D image into the clean image and two types of corruptions, namely a striped part and a laminar one.
As regularizing terms different direction dependent first and second order differences are used to cope with the specific structure of the corruptions. 
This generalizes discrete unidirectional total variational (TV) approaches. A minimizer of the model is computed by well-known primal dual techniques. 
Numerical examples show the very good performance of our new method for artificial and real-world data.
Besides FIB tomography, we have also successfully applied our technique for the removal of pure stripes in Moderate Resolution Imaging Spectroradiometer (MODIS) data.
\end{abstract}

\section{Introduction} \label{sec:intro}

The motivation for the present work was the analysis of aluminum matrix 
composites reinforced with silicon carbide particles by given high 
resolution 3D FIB tomography images.

FIB tomography, also known as serial sectioning, is an imaging technique, that has been used commercially for about $20$ years. It is applied for preparation and direct observation of structural cross-sections and the generation of microstructural data in three dimensions.
While classical $X$-ray tomography does often not reach the required material resolution, FIB tomography enables to investigate structures on a scale down to several nanometers.

\begin{figure}\centering
	\includegraphics[width=.4\textwidth]{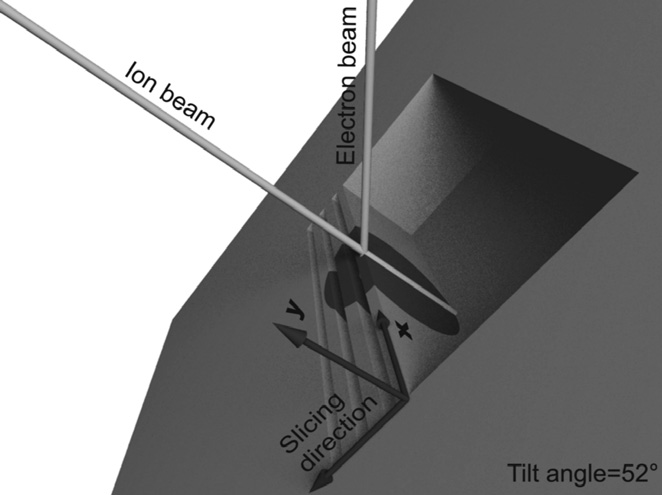}
	\caption{Illustration of a FIB system. A focused ion beam mills slice by slice of the sample in $z$-direction and the exposed surface is imaged by an electron beam. Image credit: \protect\cite{VHSSM08}}\label{fig:fib}
\end{figure}

The FIB system combines a classical scanning electron microscope with an ion beam. An illustration is shown in Fig.~\ref{fig:fib}.
The ion beam mills and polishes the material sectionally, while the electron beam is used for imaging the surface after each section.
Several hundred of these serial slices finally form a 3D image. 
For a more detailed introduction to FIB tomography we refer to \cite{BK08,GS05,Mun09}.

\begin{figure}
			\begin{tikzpicture}[every node/.style={minimum size=1cm},on grid]		
			\begin{scope}[every node/.append style={yslant=0.5,xslant=-1},yslant=0.5,xslant=-1]
			  \node at (-1.25,3.25){\includegraphics[height=2.5cm]{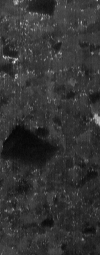}};
			\end{scope}
			\begin{scope}[every node/.append style={yslant=-0.5},yslant=-0.5]
			  \node at (-4.99,-2.99){\includegraphics[height=2.5cm]{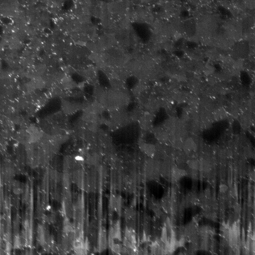}};
			  		  	\draw[thick,->,draw=white] (-3.75,-1.75) -- (-4.25,-1.75) node[text=white] at (-4.25,-2.00) {x};
			\end{scope}
			\begin{scope}[every node/.append style={yslant=0.5},yslant=0.5]
			  \node at (-3.25,0.75){\includegraphics[height=2.5cm]{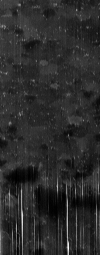}};
			  		      \draw[thick,->,draw=white] (-3.75,1.99) -- (-3.25,1.99) node[text=white] at (-3.25,1.75) {z};
			  		      \draw[thick,->,draw=white] (-3.75,1.99) -- (-3.75,1.49) node[text=white] at (-3.55,1.25) {y};
			\end{scope}
			      \node (label) at (2.3,0)[draw=none]{
			       		\includegraphics[height=.3\textwidth]{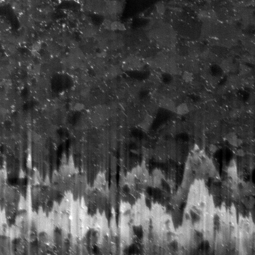} 
			       		\hspace{0.1cm}
			        	\includegraphics[height=.3\textwidth]{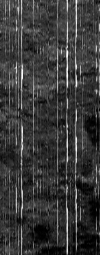} 
			        	\hspace{0.1cm}
			        	\includegraphics[height=.3\textwidth]{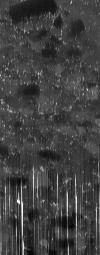}
			      };
			      \draw[thick,->,draw=white] (2.15,1.75) -- (1.65,1.75) 
			      node[text=white] at (2.15,0.95) {y};
			      \draw[thick,->,draw=white] (2.15,1.75) -- (2.15,1.25) 
			      node[text=white] at (1.35,1.75) {x};
			      \draw[thick,->,draw=white] (2.9,-1.75) -- (3.4,-1.75) 
			      node[text=white] at (3.6,-1.75) {z};
			      \draw[thick,->,draw=white] (2.9,-1.75) -- (2.9,-1.25) 
			      node[text=white] at (2.95,-1.05) {x};
			      \draw[thick,->,draw=white] (4.7,1.75) -- (5.2,1.75) 
			      node[text=white] at (5.4,1.75) {z};
			      \draw[thick,->,draw=white] (4.7,1.75) -- (4.7,1.25) 
			      node[text=white] at (4.75,0.95) {y};
			\end{tikzpicture}
			\caption{Volume image of an aluminum matrix composite obtained by FIB 
			tomography ($255 \times 255 \times 100$ pixels, i.e., $100$ slices of 
			size $11.6\mu\mathrm{m} \times 11.6\mu\mathrm{m}$) 
			and typical slices in the $x$-$y$, $x$-$z$ and $y$-$z$ 
			plane.}\label{fig:curtain}
		\end{figure}
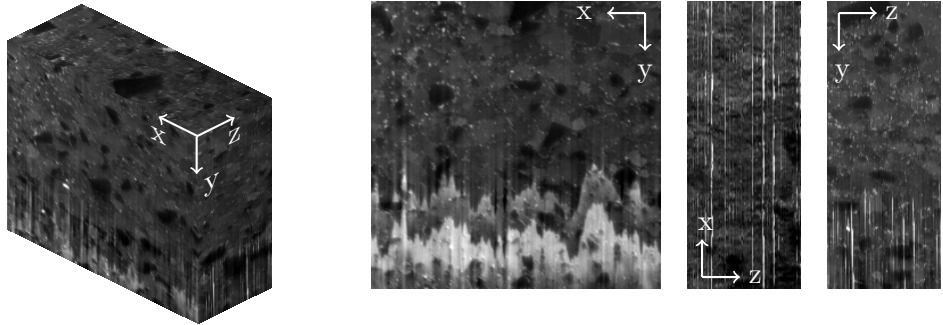

Unfortunately, such images often suffer from the so-called 
curtaining effect \cite{GS99,GS05,ZWR08}, see Fig.~\ref{fig:curtain}.
This effect arises because the sputtering rates of the ion beam are sensitive to local changes in the surface structure.
In particular below pores or cracks the rates vary and lead to stripe-like artifacts in the images.
For the aluminum matrix composite in Fig.~\ref{fig:curtain}, e.g., at the phase boundary between aluminum and the reinforced silicon carbide particles, due to their locally varying material characteristics.
Further, it may happen that the material is milled incompletely.
This leads to bright laminar artifacts.
To reduce the curtaining effects, a lower beam current could be used but this leads to a considerably higher milling time.
Note that the whole imaging and milling process already takes about $25$ hours for each data set considered here.
So further reducing the beam current leads to impractical processing times.

Therefore, we propose a computational model to reduce curtaining effects.
The corruptions in our data consist both of stripe-like structures due to different sputtering rates and laminar structures due to incomplete milling of the material.
The corrupted parts cannot be used for the further analysis of the material unless the curtaining effects are removed.
For our specimen in Fig.~\ref{fig:curtain} this was the case for more than half of the data. In FIB tomographic images the laminar corruptions are often cut off as they only occur in the lower part of the image. However, we will keep them since it is not desirable to throw away half of the data.

The problem of removing stripes in images was tackled by various approaches: 
Fourier based filtering was suggested in \cite{CSGWZ03,CP11a},
moment matching in \cite{GCS00} and histogram based methods in \cite{HW79,RTY07}.
Recently, variational models were successfully applied for the destriping of 2D images.
The central idea in \cite{BL11} is the application of unidirectional discrete TV 
terms to extract the stripes by minimizing the functional 
	\begin{align}
		\argmin_u \|\nabla_y (f-u)\|_1 + \lambda \|\nabla_{x} u\|_{1}, 
		\quad \lambda >0,
	\end{align}
	where $f$ is the original image and $u$ the destriped one.
	This method was improved in \cite{CFYL13} by adding a least squares data term and a framelet regularization term: 
	\begin{align} \label{model_4}
		\argmin_u \frac12 \|f-u\|_2^2 + \nu_1 \|\nabla_y (f-u)\|_1 &+ \nu_2 
		\|\nabla_{x} u\|_{1} + \nu_3 \|Wu\|_1,
	\end{align}
	where $W$ denotes a Parseval framelet transform.
	It has been shown that this model gives much better results than the previous one 
	also without framelet regularization, 
	i.e.~$\nu_3=0$. The framelet regularization led to another slight 
	improvement.
Finally, a more general method for the variational denoising of images with structured noise 
was developed in~\cite{FWL12}, which has also been applied for 
destriping.

Applications of the aforementioned destriping methods are, e.g., the 
restoration of
MODIS data \cite{BL11, CFYL13, RTY07}
and  FIB tomography images \cite{CFYL13,FWL12}.
In this paper, we focus mainly on the latter, 
where our FIB images are not only corrupted by 
stripe structured noise, but also by laminar effects.
To the best of our knowledge the removal of  curtaining effects involving additional large laminar parts 
was not considered before.

Compared to the above destriping methods
our new model takes the following aspects into account:
    	\begin{itemize}
    		\item[A1)]
    			The existing destriping methods remove stripes only. 
    			Since our data is not only corrupted by stripes,
    			we split the corrupted image not only into a clean part 
    			and stripes, but also into a third, laminar part.
    		\item[A2)]
    			To avoid smoothing and to keep also the small image details, 
    			we use a hard constraint for the decomposition of the corrupted image into
    			the clean, the striped and the laminar part by applying an infimal convolution model.
    		\item[A3)]
    			We consider 3D data whereas most destriping methods, except \cite{FWL12}, focus on 2D images. 
    	\end{itemize}   
	To deal with A2) we propose an infimal convolution model for the splitting.
	The infimal convolution was applied for the first 
	time in image processing by Chambolle and Lions \cite{CL97} followed by 
	many other papers, see,
	e.g., the general discrete approach in~\cite{SST11} or the continuous function space approach in \cite{HK14}.
	We also refer to a nice PhD thesis on infimal convolutions \cite{Str96}.
	Infimal convolution models were in particular used for the decomposition of images.
	
	Various variational decomposition models were proposed in the literature.
	The additive splitting into geometric and oscillating parts as texture and white Gaussian noise 
	typically exploits TV- or Besov-seminorms for the first one 
	and
	Meyers's $G$-norm \cite{AABC2003,M2001,SAC2006}, 
	the norm of the dual Sobolev space $H^{-1}$ \cite{OSV2003,VO2003}, and
	the squared $L_2$-norm of the DCT-transformed texture \cite{SMF10}
	for the second one. For specifying the noise component also the
	norm of the dual Besov space $B^\infty_{-1,\infty}$ was successfully used in \cite{AC2004}.
	Decomposition models were also applied for simultaneous structure-texture inpainting, e.g., in \cite{AGCO2006}
	as an extension of the so-called morphological component analysis model \cite{SMF10}.
	For the separation of point and curved structures based on wavelets (Besov norms) and shearlets or curvelets we refer to \cite{KL2012,SMF10}.
	{\it The clou of all these methods is the adaptation of the additive components to the task at hand.}
	For the curtaining problem we will show that our directional TV model is very well suited.
	\\
	
	{\bf Organization of the Paper.}
	In Section \ref{sec:model},
	we introduce our variational model for the removal of curtaining effects 
	and motivate its choice by the special structure of the corruptions.
	Then, in Section \ref{sec:alg}, a primal dual algorithm is presented 
	for finding a minimizer that consists of the clean image and two types of corruptions.
	Section \ref{sec:num} demonstrates the performance of our algorithm for volume images 
	obtained by FIB tomography as well as artificial examples. 
	Further, a comparison to existing methods for destriping and denoising  of 
	MODIS data is given.
	Finally, conclusions are drawn in Section \ref{sec:conc}. 

\section{New Model} \label{sec:model}
	We are given a 3D image $F \in [0,1]^{n_x,n_y,n_z}$ corrupted by curtaining effects.
	The total corruption can be seen as an overlay of two parts:
	\begin{enumerate}
		\item {\it Laminar part} consisting of very bright areas  in the $x$-$y$-plane.
		Along the $z$-axis these regions occur very fluctuant and it 
		is very unlikely that the effect appears in the same region in two subsequent 
		frames.
		\item {\it Stripes} along the $y$-axis which are few pixels wide in 
		$x$ and $z$ direction. 
	\end{enumerate}
	Based on the previous observations the corrupted image $F$ is assumed to be 
	an additive composition of the clean image $U$, the stripes $S$ and the 
	laminar part of the corruptions $L$, i.e., $F=U+S+L$.
	
	To give a sound discrete formulation of our model which takes i) and ii) into account, 
	all images are considered as column- and slice-wise reshaped vectors in $\mathbb R^N$ with $N \coloneqq 
	n_xn_yn_z$.
	The reshaped images are denoted by small letters so that in particular $f=u+s+l$.
	Then all linear operators can be written as matrices and their adjoints are just the transposed matrices.
	By $\otimes$ we denote the Kronecker product of matrices.
	We use the following directional forward difference operators:
		\begin{align}
			&{\nabla_y} \coloneqq D_y, \quad
			{\nabla_{x,z}} \coloneqq \begin{pmatrix} D_x  \\ D_z \end{pmatrix}, \quad
			{\nabla_{x,y}} \coloneqq \begin{pmatrix} D_x  \\ D_y \end{pmatrix}, \\ 
			&{\nabla_{x,y,z}} \coloneqq \begin{pmatrix} D_x  \\ D_y \\ D_z \end{pmatrix}, \quad
			{\Delta}_z \coloneqq D_{zz}, 
		\end{align}
		where
		\begin{align}
			D_x &\coloneqq I_{n_z} \otimes I_{n_y} \otimes D_{n_x}, \quad
			D_y \coloneqq I_{n_z} \otimes D_{n_y} \otimes I_{n_x}, \\
			D_z &\coloneqq D_{n_z} \otimes I_{n_y} \otimes I_{n_x}, \quad
			D_{zz} \coloneqq D^2_{n_z} \otimes I_{n_y} \otimes I_{n_x}
		\end{align}
		and
		{\small
		\begin{align}
			D_m &\coloneqq \begin{pmatrix} -1 & 1 & \\ & \ddots & \ddots \\ & & -1 & 1 \\ & & & 0 \end{pmatrix} \in \R^{m,m}, \\
			D_m^2 &\coloneqq \begin{pmatrix} 0 & 0 & \\ 1 & -2 & 1 & \\ & \ddots & \ddots & \ddots \\ & & 1 & -2 & 1 \\ & & & 0 & 0 \end{pmatrix} \in \R^{m,m}.
		\end{align}}
	Note that we assume constant extension of the image beyond the boundary for 
	the first 
	order differences and linear extension for the second order differences
	which is exactly described by the matrix form of the operators.
	Let us further emphasize that we use the Kronecker product notation only 
	for a correct description of our model 
	in a convenient matrix-vector form. In our numerical computations we work 
	with arrays based on the relation
	\begin{align}
		\vect(AXB^\tT) = (B \otimes A) \vect(X),
	\end{align}
	for $A \in \R^{m,n_x}$, $B \in \R^{k,n_y}$, $F \in \R^{n_x,n_y}$.
	
	For a vector $w \in \mathbb R^{dN}$ appearing, e.g., with $d=2$ after application of ${\nabla_{x,z}}$ or ${\nabla_{x,y}}$ and $d=3$ after applying ${\nabla_{x,y,z}}$, the grouped
	$\ell_2(\mathbb R^d)-\ell_1(\mathbb R^{N})$ norm is given by
	\begin{align} \label{norms}
		\|w\|_{2,1} \coloneqq \sum_{i=1}^{N} \sqrt{ \sum_{j=0}^{d-1} w_{i+jN}^2}.
	\end{align}
	Finally, let $\iota_\mathcal{C}$ denote the indicator function of a set $\mathcal{C}$ defined by
	$$\iota_{\mathcal{C}}(u) \coloneqq \begin{cases} 0 & \textrm{if } u \in 
	\mathcal{C}, 
	\\ +\infty & \textrm{otherwise}. \end{cases}
	$$
	We propose the following infimal convolution model
	\begin{align}\label{IC-formulation}
		\inf_{u+s+l=f} {\varphi_1(u) + \varphi_2(s) + 
		\varphi_3(l)},
	\end{align}
	with
	\begin{align}
		\varphi_1(u) &\coloneqq \mu_1 \|\nabla_{x,z} u\|_{2,1} + \mu_2 \|\Delta_z 
		u\|_1 
		+ \iota_{[0,1]^N}(u),\\
		\varphi_2(s) &\coloneqq \|\nabla_y s\|_1,\\
		\varphi_3(l) &\coloneqq \mu_3 \|\nabla_{x,y} l\|_{2,1}	
	\end{align}
	and regularization parameters $\mu_1,\mu_2,\mu_3 >0$.
	For convenience, we also add the definition of the above norms for the corresponding array formulation:%
	{\small
	\begin{align}
		\|\nabla_{x,z} U \|_{2,1} &\coloneqq \sum_{i,j,k=1}^{n_x,n_y,n_z} \Big( \left(U_{i+1,j,k}-U_{i,j,k} \right)^2 + (U_{i,j,k+1}-U_{i,j,k})^2\Big)^{\frac12},\\
		 \| \Delta_{z} U \|_{1} &\coloneqq\sum_{i,j,k=1}^{n_x,n_y,n_z} |U_{i,j,k+1}-2 U_{i,j,k}+ U_{i,j,k-1}|,\\ 
		 \| \nabla_{y} S \|_{1} &\coloneqq \sum_{i,j,k=1}^{n_x,n_y,n_z} |U_{i,j+1,k}-U_{i,j,k}|%
	\end{align}}%
	with the appropriate mirrored or linear extensions at the boundary.
	The last summand in $\varphi_1$ ensures that the clean image $u$ has its range in $[0,1]$.
	The choice of the other terms is motivated by i) and ii) as follows:
	\begin{itemize}
	 \item[i)] $\|\nabla_{x,y} l\|_{2,1}$ and $\|\Delta_z u\|_1$: {\it laminar 
	 distortions} have a small width in 
	$z$ direction, thus a high response in the second order derivative is 
	expected along this direction. Therefore, we penalize second order 
	differences in $z$ direction of $u$.  Additionally, these corruptions are 
	spacious in $x$ and $y$ direction such 
	that a bidirectional TV regularizer is useful for $l$.
	\item[ii)]  $\|\nabla_y s\|_1$ and $\|\nabla_{x,z} u\|_{2,1}$:  {\it stripes} 
	 occur only in $y$ direction so that we penalize first order differences 
	in $y$ direction of $s$. This was also done in \cite{BL11,CFYL13} setting 
	$s = f-u$.
	Further, we penalize coupled first order differences in $x$ and $z$ 
	direction of 
	$u$ since we assume that the main differences along these directions occur 
	in $s$ while $u$ has only few edges.
	\end{itemize}
	
	\begin{remark} \label{rem:TV}
	        Alternatively to $\mu_1 \| \nabla_{x,z} u \|_{2,1} + \mu_2 \|\Delta_z u\|_1$ 
		we can also use the 3D TV term $\| \nabla_{x,y,z} u \|_{2,1}$, i.e., instead of  \eqref{IC-formulation}
		we are asking for
		\begin{align}\label{IC-formulation_rev}
				\inf_{u+s+l=f} {\varphi_1(u) + \varphi_2(s) + \varphi_3(l)},
		\end{align}
		where $\varphi_i$, $i=2,3$ are defined as in \eqref{IC-formulation} and
		\begin{align}
			\varphi_1(u) &\coloneqq \mu_1 \|\nabla_{x,y,z} u\|_{2,1}
			+ \iota_{[0,1]^N}(u).
		\end{align}
		The term for $u$ corresponds to the assumption that $u$ is a piecewise smooth image.
		Our model \eqref{IC-formulation} focuses more on the corruptions. 
		We argue in i) that, in contrast to the laminar corruptions, $u$ has a small second order derivative in $z$ direction. 
		This viewpoint is also taken in \cite{BL11,CFYL13}. We have implemented both methods.		 
		The difference between the models is best visible in Fig.~\ref{fig:toy} in the numerical section.
	\end{remark}
	
	Proposition~\ref{exist} shows that our models \eqref{IC-formulation} and \eqref{IC-formulation_rev} have a solution.
	Its uniqueness cannot be expected since arbitrary constants can be added 
	to $s$ and subtracted from $l$ without changing the objective function.

	\begin{proposition}\label{exist}
		Let $f \in [0,1]^{N}$. Then  \eqref{IC-formulation}, respectively \eqref{IC-formulation_rev}, has  a minimizer.
	\end{proposition}
	The proof is given in the appendix.

\section{Algorithm} \label{sec:alg}
	In this section, we propose to use a primal dual algorithm to solve the convex, but non-smooth problem \eqref{IC-formulation}.
	Problem \eqref{IC-formulation_rev} can be handled in a similar way.
	Primal dual algorithms have recently found wide applications in image processing and machine learning. 
	The reason is that these algorithms split the original task  into a sequence of simple computable, proximal mappings.
	This paper cannot address the huge amount of papers on various primal dual strategies, which were successfully applied
	in imaging. For a good overview on primal dual methods we refer to \cite{CP11b,KP14}. 
	
	Recall that the proximal mapping of a proper, convex, lower semi continuous
	function $\varphi \colon \mathbb R^N \rightarrow (-\infty, + \infty]$ is defined by
	\begin{align} \label{eq:prob}
	\prox_{\lambda \varphi} (x) \coloneqq \argmin_{y \in \R^N} \left\{ \varphi (y) + \frac{1}{2\lambda} \|x-y\|_2^2 \right\}.
	\end{align} 
	Indeed the proximal mapping at $x \in \mathbb R^N$ exists and is unique.
	To make primal dual algorithms efficient, the proximal mapping should in addition be simply computable.
	For applications of proximal mappings see, e.g., \cite{PB2013}.
	
	We want to apply the primal dual algorithm with an extrapolation of the dual variable 
	suggested by Chambolle, Pock and Cremers in \cite{CP11,PCCB09}, see also \cite{CP15,GEB13,PC11} for extensions.
	We call this algorithm PDHG (Primal Dual Hybrid Gradient) algorithm in accordance with the name of a similar algorithm which was presented
	without the extrapolation idea and without convergence proof in \cite{ZC08}. Later its convergence was examined in \cite{BR2012}.
	
	To this end, we have to rewrite our model \eqref{IC-formulation}.
	Using the vectors $x \coloneqq (u,s,l)^\tT$ and $y \coloneqq (y_i)_{i=1}^4$ and the notation
	\begin{align}\label{E}
	\mathcal{C} &\coloneqq \{  (u,s,l) : u+s+l=f, \; u \in [0,1] \},\\
	h(y) &\coloneqq \mu_1 \| y_1 \|_{2,1} + \mu_2 \| y_2 \|_1 +  \| y_3 \|_{2,1} + \mu_3 \|y_4\|_1
	\end{align}
	we ask for 
	\begin{align}
	\argmin_{x,y} \left\{ \iota_{\mathcal{C}}(x) + h(y) \right\} \quad \mathrm{subject \; to} \quad K x =y,
	\end{align}
	where
	\begin{align}
		 K \coloneqq	\begin{pmatrix}
			\nabla_{x,z} & 0 & 0 \\
			\Delta_z & 0 & 0\\
			0 & \nabla_y & 0 \\
			0 & 0 & \nabla_{x,y} 
			\end{pmatrix}.			
	\end{align}
	Then the PDHG algorithm from \cite{CP11,PCCB09} with an extrapolation of the dual variable $p$ reads
	\begin{align}
		x^{k+1} &= \prox_{\tau \iota_\mathcal{C}} (x^{k} - \tau K^\tT \bar p^{k}) = \mathcal{P}_\mathcal{C} (x^{k} - \tau K^\tT \bar p^{k}),\\
		p^{k+1} &= \prox_{\sigma h^*} (p^{k} + \sigma K x^{k+1}),\\
		\bar p^{k+1} &= p^{k+1} + \theta(p^{k+1} - p^{k}), \quad \theta \in (0,1],
	\end{align}
	where $h^*$ denotes the conjugate function of $h$ and $\mathcal{P}_\mathcal{C}$ is the orthogonal projection onto $\mathcal{C}$.
	The algorithm converges if the parameters fulfill $\tau \sigma < \frac{1}{\| K \|_2^2}$.
	For our structured matrix $K$ it is not hard to check by diagonalizing it with trigonometric transforms
	that $\|K\|_2^2 \le 24$, so that $\tau \sigma < \frac{1}{24}$ must be fulfilled.
	
	Applying Moreau's identity, one can directly work with the function $h$ instead of its conjugate. 
	Using the rescaling $b^{k} \coloneqq p^{k}/\sigma$, the algorithm becomes
	\begin{align}
		x^{k+1} &= \mathcal{P}_\mathcal{C} (x^{k} - \tau K^\tT \bar p^{k}),\\
		y^{k+1} &= \prox_{\frac{1}{\sigma} h} (b^{k} + K x^{k+1}),\\
		b^{k+1} &= b^{k} + K x^{k+1} - y^{k+1},\\
		\bar b^{k+1} &= b^{k+1} + \theta (b^{k+1} - b^{k}), \quad \theta \in (0,1].
	\end{align} 
	For details, see, e.g., \cite{BSS14}. Fortunately, the proximal mapping of $h$ can be separated into the proximal mappings
	of its summands with respect to $y_i$, $i=1,\ldots,4$. The whole algorithm with our parameter choice is detailed in Algorithm \ref{algo}.

	{\small
	\begin{algorithm} \label{algo}
		\KwInit{$u^{(0)}=f$, $s^{(0)}=0$, $l^{(0)}=0$, 
		$b_i^{(0)}=0$, $\bar b_i^{(0)}=0$, $i=1,2,3,4$, $\theta =1$, $\tau = 
		\frac15$, $\sigma=\frac15$.} \\
		\KwIter{For $k = 0,1,\ldots$ iterate}
		\begin{align}
			(u,s,l)^{(k+1)} = \ &\mathcal{P}_{\mathcal{C}}\big(
			(u^{(k)}-\tau \sigma {\nabla}_{x,z}^\tT \bar b_1^{(k)}-\tau \sigma 
			{\Delta}_{z}^{\tT} \bar b_4^{(k)},\\
			& \quad s^{(k)}-\tau \sigma {\nabla}_y^\tT \bar 
			b_2^{(k)}, l^{(k)}-\tau \sigma {\nabla}_{x,y}^\tT \bar 
			b_3^{(k)})\big)
			\\
			y_1^{(k+1)} = \ &\prox_{\frac{\mu_1}{\sigma} \|\cdot\|_{2,1}} 
			(b_1^{(k)} + {\nabla}_{x,z} u^{(k+1)})\\
			y_2^{(k+1)} = \ &\prox_{\frac{\mu_2}{\sigma} \|\cdot\|_{1}} 
			(b_2^{(k)} + {\Delta}_{z} u^{(k+1)})
			\\
			y_3^{(k+1)} = \ &\prox_{\frac1\sigma \|\cdot\|_{1}}
			(b_3^{(k)} + {\nabla}_y s^{(k+1)})
			\\
			y_4^{(k+1)} = \ &\prox_{\frac{\mu_3}{\sigma} \|\cdot\|_{2,1}} 
			(b_4^{(k)} + {\nabla}_{x,y} l^{(k+1)})
			\\
			b_1^{(k+1)} = \ & b_1^{(k)} + {\nabla}_{x,z}   u^{(k+1)} - 
			y_1^{(k+1)}	
			\\
			b_2^{(k+1)} = \ & b_2^{(k)} + {\Delta}_z   u^{(k+1)} - 
			y_2^{(k+1)}
			\\
			b_3^{(k+1)} = \ & b_3^{(k)} + {\nabla}_{x,y}   l^{(k+1)} 
			- y_3^{(k+1)}
			\\
			b_4^{(k+1)} = \ & b_4^{(k)} + {\nabla}_y   s^{(k+1)} - y_4^{(k+1)}
			\\ 
			\bar b_i^{(k+1)} = \ & b_i^{(k+1)} + \theta (b_i^{(k+1)}-b_i^{(k)})	
			\quad \ \ i=1,\ldots,4		
		\end{align}
		\KwOutput{{\rm Clean image} $u$, {\rm corruptions} $s$, $l$}
		\caption{PDHG for  \protect\eqref{E}.} \label{alg1}
	\end{algorithm}
	}

	The first step of the algorithm is a voxelwise projection
	$\mathcal{P}_\mathcal{C}$ onto the set
	$\mathcal{C}$. Since $\mathcal{C}$ is just a plane with an additional range constraint for 
	$u$, this projection can be computed in a straightforward way.
	The update steps for $y_2$ and $y_3$ require a componentwise soft shrinkage, which is defined for $t \in \R$ by
	\begin{align}
		S_\lambda(t) \coloneqq &
			\begin{cases}
					0 & \mathrm{if} \ |t| \le \lambda,\\
					t(1- \tfrac{\lambda}{|t|}) & \mathrm{if} \ |t| > \lambda.
			\end{cases}
	\end{align}
	For $y_1$ and $y_4$ a componentwise coupled shrinkage procedure is required. 
	For vectors $\mathbf{t} = (t_1,t_2)^\tT \in \R^2$ it is given by
	\begin{align}
		\mathbf{S}_\lambda(\mathbf{t}) \coloneqq \begin{cases}
						0 & \mathrm{if} \ \|\mathbf{t}\|_2 \le \lambda,\\ 
						\mathbf{t}(1-\tfrac{\lambda}{\|\mathbf{t}\|_2}) & \mathrm{if} \ \|\mathbf{t}\|_2 > \lambda.
		\end{cases}
	\end{align}
	
	Due to the separation of the computation of the $y_i$, $i=1,\ldots,4$,
	the algorithm can be implemented in parallel on a multi-core architecture.
	Recently such  parallel strategies were also examined by 
	Pesquet and co-authors \cite{BCPP11,PCP11}.
	Further we want to mention the interesting stochastic block coordinate 
	algorithms in \cite{CP2015,RCP2015}.	
	
\section{Numerical Examples} \label{sec:num}
In this section, we present the results of our algorithm, which was implemented in MATLAB2015b,
for 3D artificial data,  3D FIB tomography data 
and 2D MODIS data. 
For the 3D examples, videos of the whole volume images are  provided at our website {\small \url{http://www.mathematik.uni-kl.de/imagepro/members/fitschen/curtaining0/}}.

We compare the results of our model \eqref{IC-formulation} with three other methods M1, M2 and M3
and with the  modified model \eqref{IC-formulation_rev} in Remark \ref{rem:TV}.
Since there seems to be no work especially on curtaining effects, the proposed 
methods are compared to the following 3D destriping and denoising algorithms:
\begin{enumerate}
	\item[M1)]	Firstly, we generalize the 2D destriping approach from \cite{CFYL13} with regularization terms as 
		in \eqref{model_4} without the framelet part to 3D data, which results in
		\begin{align}
			\argmin_u \frac12 \|f-u\|_2^2 + \nu_1 \|\nabla_y (f-u)\|_1 + \nu_2 
			\|\nabla_{x,z} u\|_{2,1}.
		\end{align}
		To remove the corruptions that are not stripes, a median filter in 
		$z$ direction is applied beforehand since the simple destriping 
		algorithm is not able to remove these artifacts.
		
	\item[M2)]   Secondly, we consider the 2D  denoising method proposed in \cite{FWL12} for removing structured 
		noise. We have used the implementation available at 
		\url{http://www.math.univ-toulouse.fr/~weiss/PageCodes.html}\\ with a 
		Gabor filter and parameters $(p,\alpha)$.
		Notice that this method is applicable for a wide range of structured noise. Thus 
		one can not expect that the results are as good those of our model which is 
		adjusted to this specific task.
		Again, we apply a median filter in $z$ direction first and then apply 
		the algorithm to the 2D slices 
		separately.
	\item[M3)] Thirdly, we consider the 3D shearlet hard thresholding method in \cite{KLR14} for denoising of videos, see also \cite{HSHPSLS12}. We have used the software package "ShearLab 3D" available at \url{http://www.shearlab.org/software}.
	For this method we have used the shearlet system $SL3D_1$, i.e., $3$ scales and a redundancy of $76$.
	In total we have used $5$ different thresholding factors, one factor $v_1$ for the direction orthogonal to the laminar part, another one $v_2$ for the directions nearest to that direction, in the same way another two for the stripes ($v_3, v_4$) and finally a fifth factor $v_5$ for all other directions.
\end{enumerate}		
Compared to the methods M1 and M2 our method has the advantage of computing everything in one step without any preprocessing. 

\subsection*{Artificial Data}
The first two experiments deal with artificial 3D images disturbed by 
stripes and laminar structures, that are fluctuant in $z$ direction.
The parameters of all models were chosen by an extensive grid search such that the PSNRs were optimized.

\begin{figure*}	\centering \begin{tabular}{ccc}
		corrupted & proposed \eqref{IC-formulation}& proposed \eqref{IC-formulation_rev} \\ 
		\includegraphics[width=.3156\textwidth]{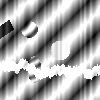}
		&\includegraphics[width=.3156\textwidth]{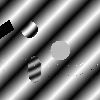}
		&\includegraphics[width=.3156\textwidth]{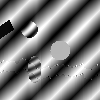}\\
		 M1 & M2 & M3\\
		\includegraphics[width=.3156\textwidth]{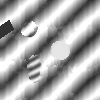}
		&\includegraphics[width=.3156\textwidth]{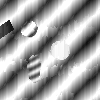}
		&\includegraphics[width=.3156\textwidth]{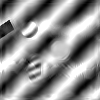}
		\end{tabular}
		\caption{Results of various methods for the artificial image. 
		Top left to bottom right: 
		Exemplary $x$-$y$ slice of the corrupted image, reconstructed images 
		with our method \protect\eqref{IC-formulation} ($(\mu_1,\mu_2,\mu_3) = 
		(\frac{1}{1500},\frac{4}{300},\frac{7}{300})$, PSNR $34.62$), 
		our method \protect\eqref{IC-formulation_rev}  ($(\mu_1,\mu_3) = 
		(\frac{1}{300},\frac{3}{300})$, PSNR 27.74),
		the generalized destriping method M1 ($(\nu_1,\nu_2) = (0.9, 
		0.5)$, PSNR 14.31), the generalized 
		denoising algorithm M2 ($p=1, \alpha = 0.13$, PSNR 
		14.05) and the shearlet hard thresholding method ($v_1=5, v_2=3, v_3=1.1, v_4=1, v_5=\frac13$).}\label{fig:toycomp}
\end{figure*}

In Fig.~\ref{fig:toycomp}, a comparison of the 
proposed models \eqref{IC-formulation} and \eqref{IC-formulation_rev} and the methods M1, M2 and M3 is shown.
Specific slices of the volume image are presented. 
The method \eqref{IC-formulation} performs clearly better than 
the other ones. 
This is also confirmed by the corresponding PSNR values for the whole 3D data.
If we apply M1 and M2 without a median filter in the preprocessing 
step, the PSNR is even worse, namely PSNR = 12.53 for the destriping method and PSNR = 10.05 for the denoising method.
In contrast to these methods with M3 it is possible to remove the laminar part without an additional median filter. 
But nevertheless the results are clearly worse compared to the proposed model.

\begin{figure*} \centering \begin{tabular}{cccc}
		clean & corrupted & proposed \eqref{IC-formulation}& proposed \eqref{IC-formulation_rev} \\
		\includegraphics[width=.23\textwidth]{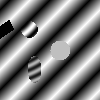}
		&\includegraphics[width=.23\textwidth]{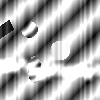}
		&\includegraphics[width=.23\textwidth]{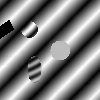}
		&\includegraphics[width=.23\textwidth]{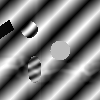}\\[1ex]
		\includegraphics[width=.23\textwidth]{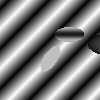}
		&\includegraphics[width=.23\textwidth]{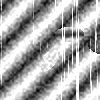}
		&\includegraphics[width=.23\textwidth]{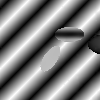}
		&\includegraphics[width=.23\textwidth]{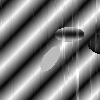}\\[1ex]
		\includegraphics[width=.23\textwidth]{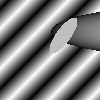}
		&\includegraphics[width=.23\textwidth]{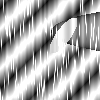}
		&\includegraphics[width=.23\textwidth]{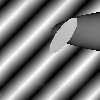}
		&\includegraphics[width=.23\textwidth]{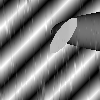}
		\end{tabular}
		\caption{Results of the proposed methods for an artificial image. Top to bottom: exemplary slices in $x$-$y$, $x$-$z$ and $y$-$z$ 
		directions of the volume image. 
		Left to right: original image, corrupted image, reconstructed image $u$ with our method \protect\eqref{IC-formulation} ($(\mu_1,\mu_2,\mu_3) = 
				(\frac{1}{1500},\frac{4}{300},\frac{7}{300})$, PSNR 32.76) and \protect\eqref{IC-formulation_rev}  ($(\mu_1,\mu_3) = (\frac{1}{500},\frac{2}{300})$, PSNR 27.39).}\label{fig:toy}
\end{figure*}

\setlength{\tabcolsep}{.1cm}
Comparing the models \eqref{IC-formulation} and \eqref{IC-formulation_rev}, we observe that the main differences occur along the laminar corruptions.
To exclude that the reason for this behavior are the sharp edges of the laminar part, we  show another example with a smoothed laminar part in
Fig.~\ref{fig:toy}.
Model \eqref{IC-formulation} gives an almost perfect reconstruction, whereas the method \eqref{IC-formulation_rev} is still not able to remove the laminar corruptions correctly.
The results are summed up in Table \ref{tab:error}.
\begin{table}
\begin{tabular}{cc|ccccc}
	Image & Error measure & Model \eqref{IC-formulation} & Model \eqref{IC-formulation_rev} & M1 & M2 & M3 \\ \hline
							& PSNR & $34.62$ & $27.74$ & $14.31$ & $14.05$ & $12.17$ \\
	Fig.~\ref{fig:toycomp} 	 & MSE & $0.35\cdot 10^{-3}$ & $1.68\cdot 10^{-3}$ & $0.0371$ & $0.0394$ & $0.0607$\\
						 	 & SSIM & $0.99$ & $0.959$ & $0.805$ & $0.784$ & $0.741$\\ \hline
							& PSNR & $32.76$ & $27.39$ & $14.14$ & $13.92$ & $12.22$\\
	Fig.~\ref{fig:toy} 		 & MSE & $0.53\cdot 10^{-3}$ & $1.82\cdot 10^{-3}$ & $0.0386$ & $0.0406$ & $0.06$\\
							 & SSIM & $0.99$ & $0.959$ & $0.808$ & $0.782$ & $0.7512$
\end{tabular}
\caption{Error measures for the results of various methods.}\label{tab:error}
\end{table}
Besides the PSNR, which we have discussed, the mean square error (MSE) and structural similarity index (SSIM) are stated in the table.
In all cases our method \eqref{IC-formulation} performs best and the methods M1, M2 and M3 are clearly worse.

\subsection*{FIB Data}
Next we apply the methods to different FIB tomography data, namely
to aluminum matrix composites with different particle sizes
and to a bronze sample.
The images of the aluminum matrix composite were obtained from 
the ``Nano Structuring Center'' Kaiserslautern
and the bronze sample from the ``Material Engineering Center Saarland'' in Saarbr\"ucken.

For all data sets we applied the proposed model \eqref{IC-formulation} with the same parameters, namely
\begin{align}
	(\mu_1,\mu_2,\mu_3) = (\tfrac{1}{300},\tfrac{2}{300},\tfrac{6}{300}).
\end{align}
Only for the bronze sample we set $\mu_2 = 0$ since there is no laminar part.

\begin{figure*}\centering\begin{tabular}{cccc}
		corrupted & reconstructed & stripes & laminar\\
		\includegraphics[width=.23\textwidth]{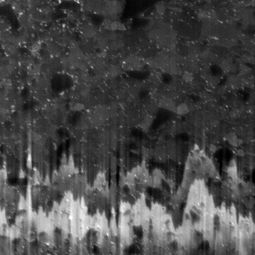}
		&\includegraphics[width=.23\textwidth]{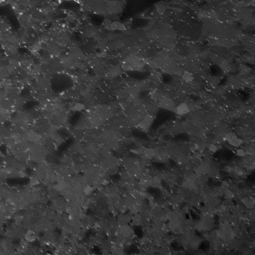}
		&\includegraphics[width=.23\textwidth]{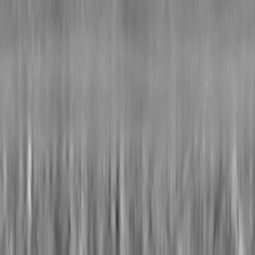}
		&\includegraphics[width=.23\textwidth]{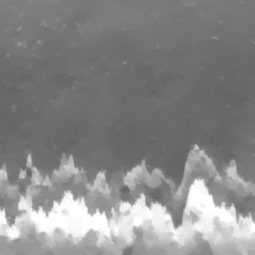}\\[1ex]
		\includegraphics[width=.23\textwidth]{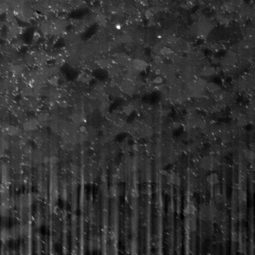}
		&\includegraphics[width=.23\textwidth]{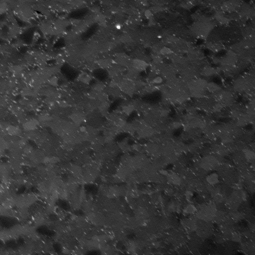}
		&\includegraphics[width=.23\textwidth]{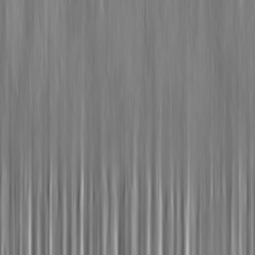}
		&\includegraphics[width=.23\textwidth]{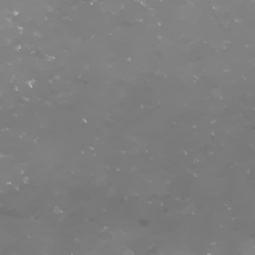}	
		\end{tabular}
		\caption{Results of our method \protect\eqref{IC-formulation} for real FIB data. 
		Left to right: Exemplary $x$-$y$ slices of the corrupted image 
		$F$, reconstructed image $U$, corruptions $S$ and $L$. The first row shows a slice 
		with prominent corruptions also in $L$ and the second row with almost only stripes.}\label{fig:real}
\end{figure*}

Fig.~\ref{fig:real} demonstrates that our method \eqref{IC-formulation} 
is able to split the corrupted image nicely into a clean image and two corrupted parts.

\begin{figure*}\centering\begin{tabular}{ccccc}
		original & proposed \eqref{IC-formulation}& proposed \eqref{IC-formulation_rev} & M1 & M3 \\
		\includegraphics[width=.185\textwidth]{./images/real_ex/real_xy_corr.png}
		&\includegraphics[width=.185\textwidth]{./images/real_ex/real_xy_rec.png}
		&\includegraphics[width=.185\textwidth]{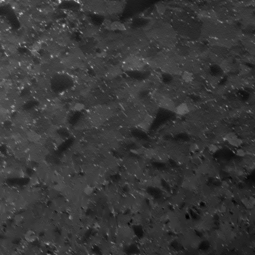}
		&\includegraphics[width=.185\textwidth]{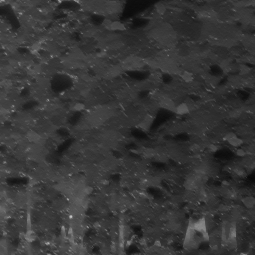}
		&\includegraphics[width=.185\textwidth]{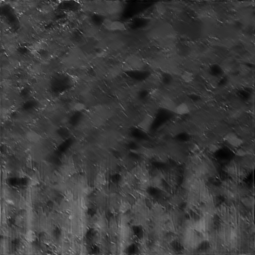}\\[1ex]
		\includegraphics[width=.185\textwidth]{./images/real_ex/real_xy_corr_s.png}
		&\includegraphics[width=.185\textwidth]{./images/real_ex/real_xy_rec_s.png}
		&\includegraphics[width=.185\textwidth]{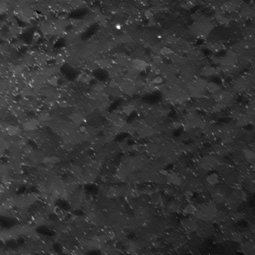}
		&\includegraphics[width=.185\textwidth]{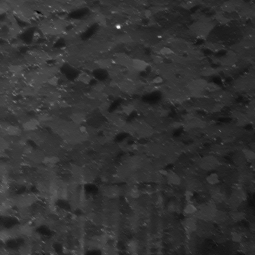}
		&\includegraphics[width=.185\textwidth]{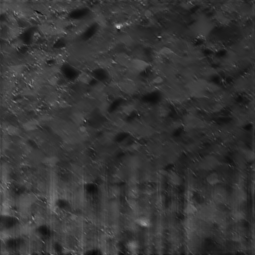}
		\end{tabular}
		\caption{Results of various methods for real FIB data.
		Left to right: Exemplary $x$-$y$ slices of the corrupted image, 
		reconstructed image with the proposed model \protect\eqref{IC-formulation}, the proposed model \protect\eqref{IC-formulation_rev} ($(\mu_1,\mu_3) = (\frac{3}{300},\frac{6}{300})$), the destriping method M1 ($(\nu_1,\nu_2) = (0.4, 0.2)$) and shearlet thresholding M3 ($v_1=4, v_2=3, v_3=3, v_4=3, v_5=\frac{1}{30}$).}\label{fig:real2}
\end{figure*}

In Fig.~\ref{fig:real2}, a comparison of our model  \eqref{IC-formulation}  with \eqref{IC-formulation_rev}, M1 and M3 is
presented. 
For the slice that is mainly corrupted by stripes, M1 shows acceptable results, 
but cannot compete with \eqref{IC-formulation} and \eqref{IC-formulation_rev}, which perform equally well here.
If there are not only stripes, one has to apply a median filter beforehand to obtain reasonable results with M1. 
However, this results in a loss of small image structures
like the bright aggradations in the upper left part of the images.
The method M3 leads to slightly smoothed results and still shows some remains of both stripes and the laminar part.

\begin{figure*}\centering \begin{tabular}{cccc}
		original & proposed \eqref{IC-formulation} & M1 & M3 \\
		\includegraphics[width=.23\textwidth]{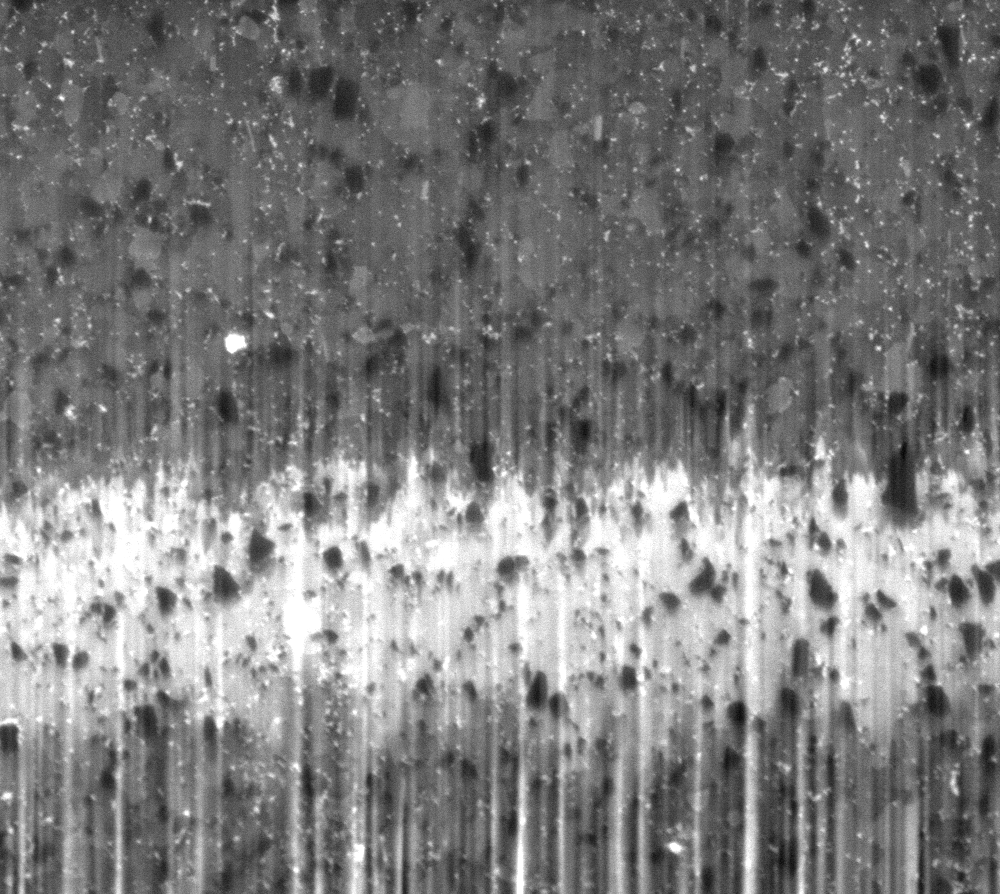}
		&\includegraphics[width=.23\textwidth]{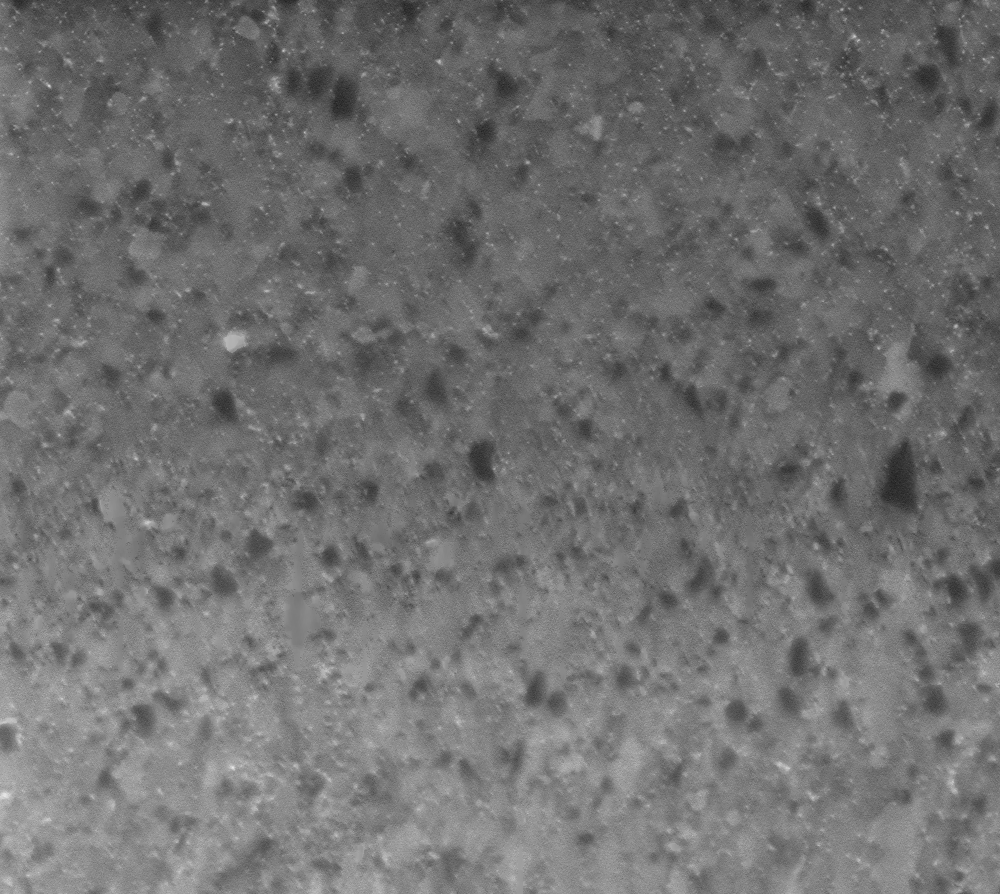}
		&\includegraphics[width=.23\textwidth]{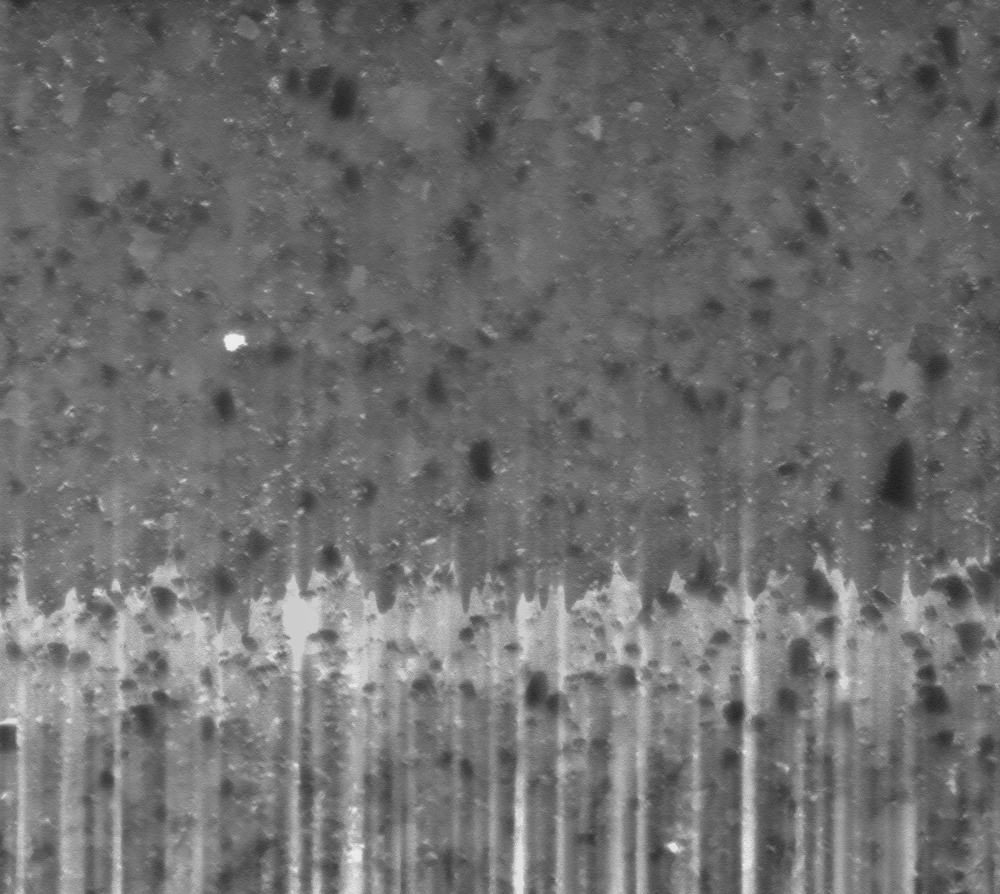}
		&\includegraphics[width=.23\textwidth]{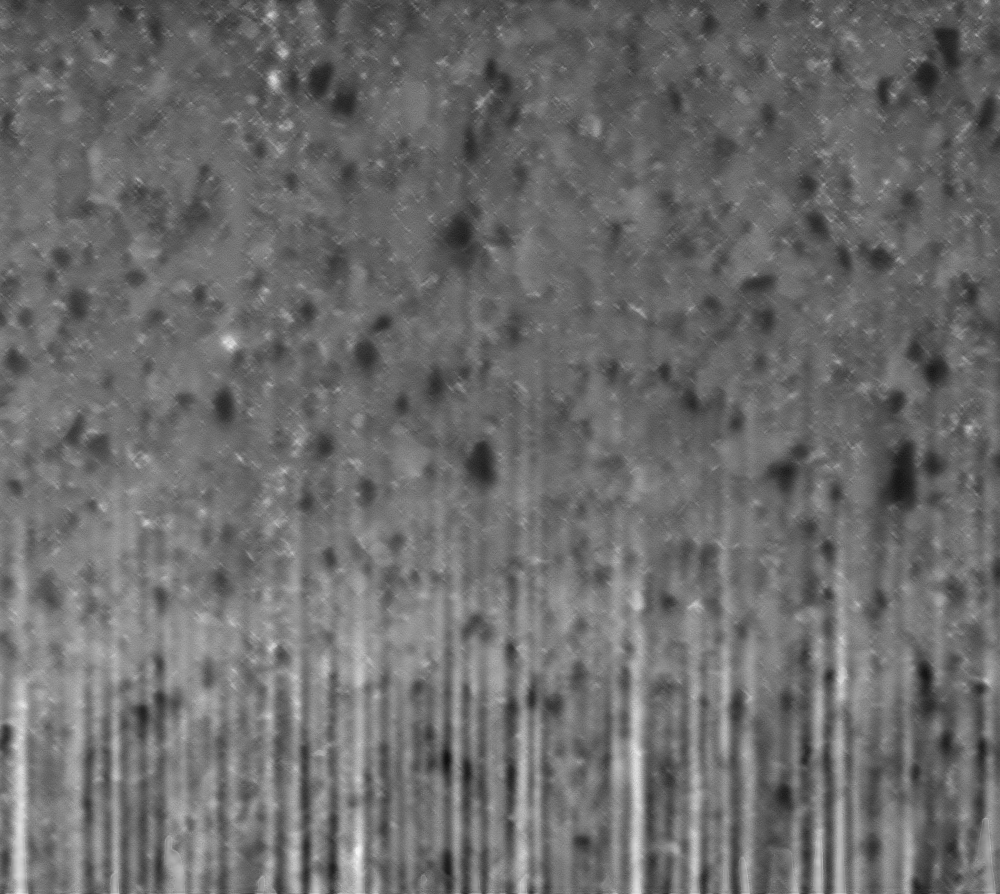}\\[1ex]
		\includegraphics[width=.23\textwidth]{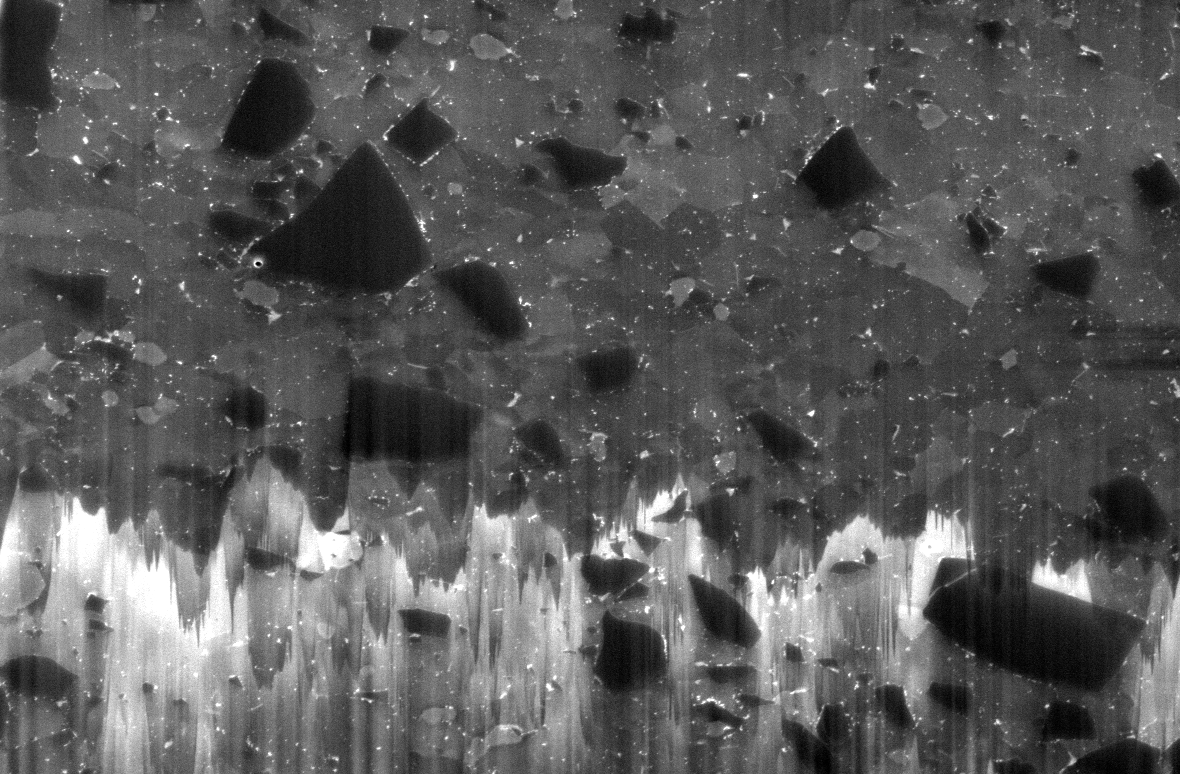}
		&\includegraphics[width=.23\textwidth]{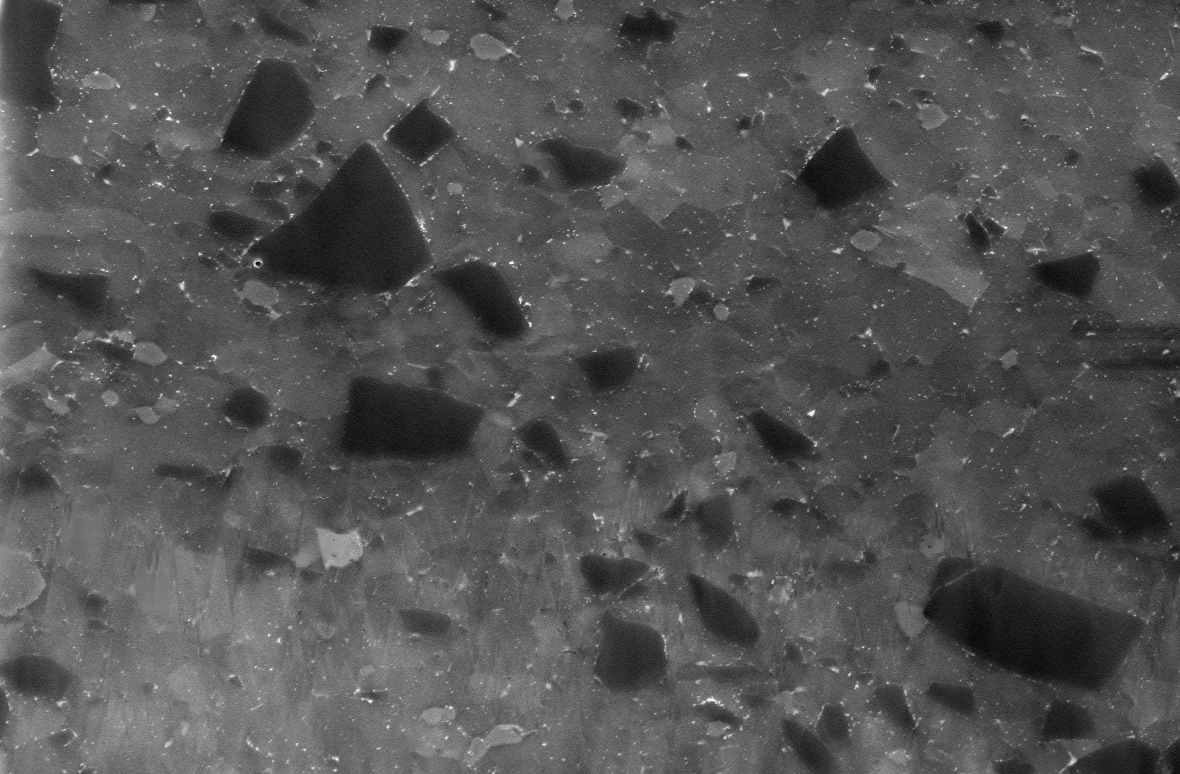}
		&\includegraphics[width=.23\textwidth]{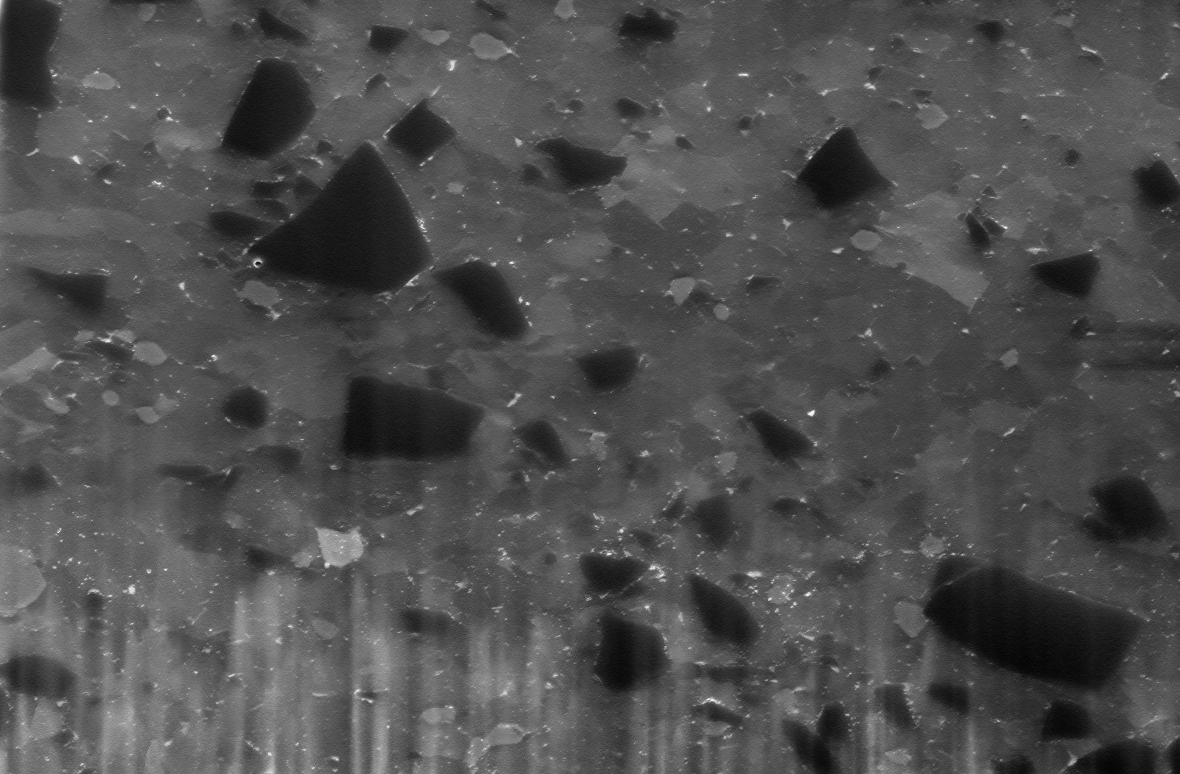}
		&\includegraphics[width=.23\textwidth]{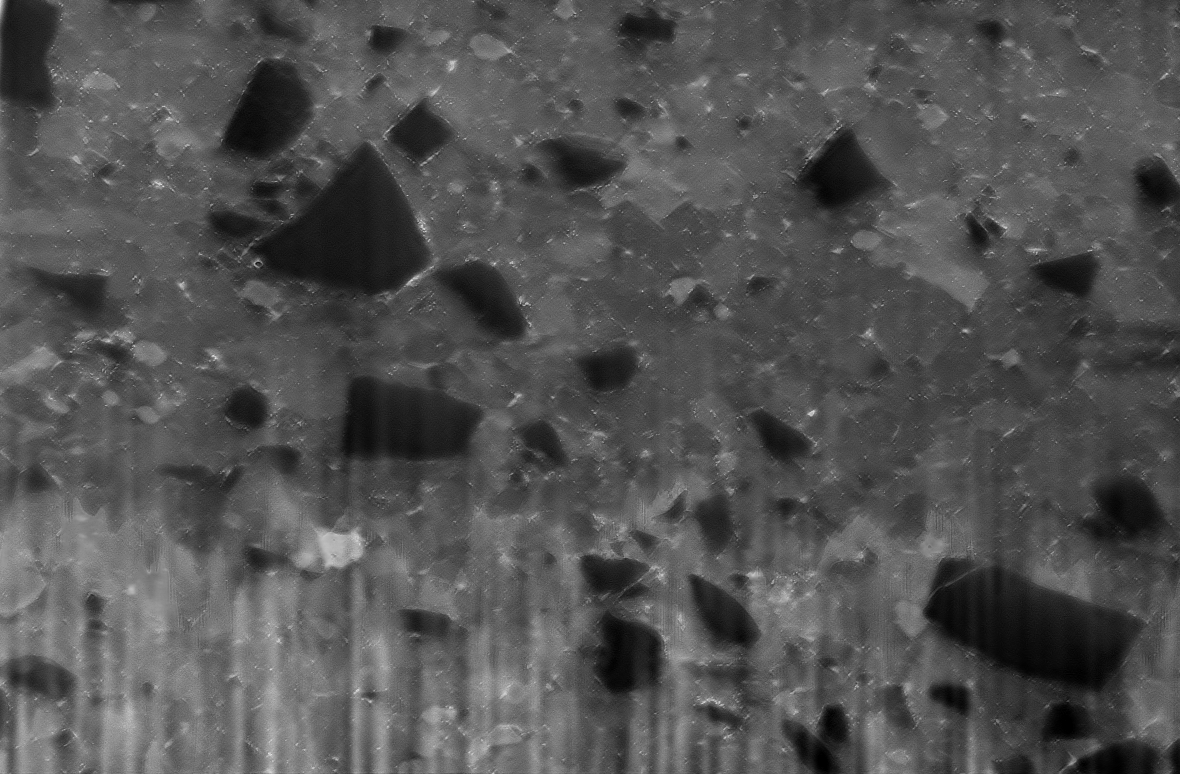}\\[1ex]
		\includegraphics[width=.23\textwidth]{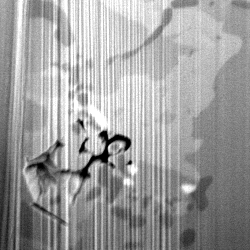}
		&\includegraphics[width=.23\textwidth]{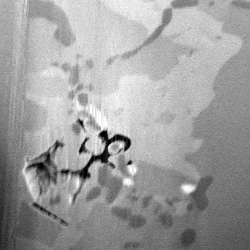}
		&\includegraphics[width=.23\textwidth]{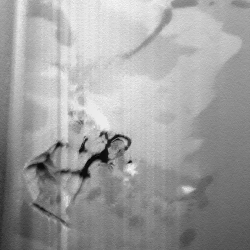}
		&\includegraphics[width=.23\textwidth]{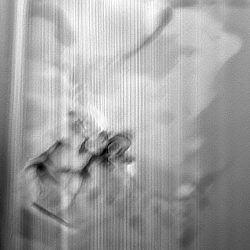}
		\end{tabular}
		\caption{Results of our method \protect\eqref{IC-formulation} and M1 ($(\nu_1,\nu_2) = (0.4, 0.2)$) for different materials. From top to bottom: aluminum matrix composite with small particles, large particles and a bronze sample.}\label{fig:real3}
\end{figure*}

Finally, we show results obtained by our model \eqref{IC-formulation}
for FIB images of various materials. 
In Fig.~\ref{fig:real3} exemplary slices are shown. 
We obtain very good results although the materials as well as the corruptions look different.
For a comparison also the results of M1 and M3 are presented. However, these are clearly worse. In particular with M1 the laminar structures are still present and M3 is not able to remove the stripes. But note that both methods are not designed for such corruptions.

\subsection*{MODIS Data}
\begin{figure*}\centering\begin{tabular}{cccc}
		original & proposed \eqref{IC-formulation}\\
		\includegraphics[width=.23\textwidth, angle =90]{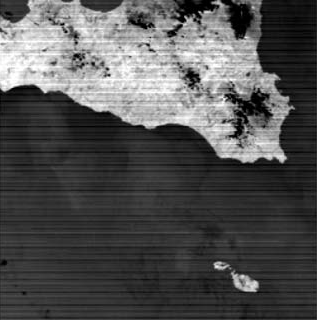}
				\includegraphics[width=.23\textwidth, angle =90]{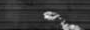}\hspace{0.1cm}
		&\includegraphics[width=.23\textwidth, angle =90]{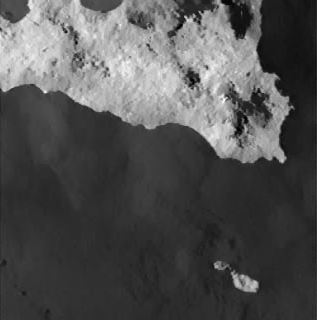}
				\includegraphics[width=.23\textwidth, angle =90]{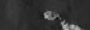}\\[1.5ex]
				 M1 & M2 \\
		\includegraphics[width=.23\textwidth, angle =90]{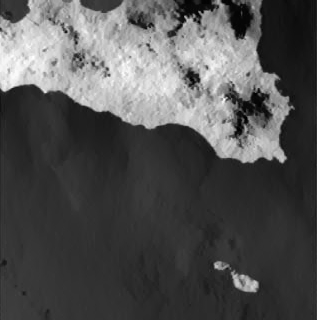}
				\includegraphics[width=.23\textwidth, angle =90]{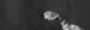}\hspace{0.1cm}
		&\includegraphics[width=.23\textwidth, angle =90]{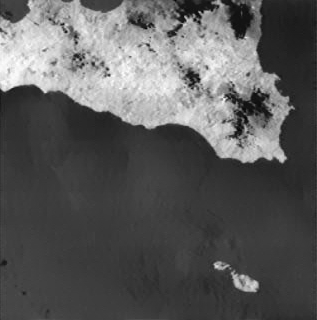}
				\includegraphics[width=.23\textwidth, angle =90]{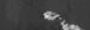}\\[1.5ex]
		original & proposed \eqref{IC-formulation} \\
		\includegraphics[width=.23\textwidth, angle =90]{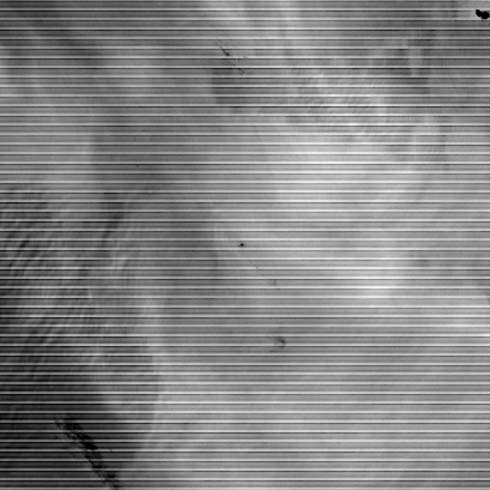}
				\includegraphics[width=.23\textwidth, angle =90]{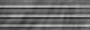}\hspace{0.1cm}
		&\includegraphics[width=.23\textwidth, angle =90]{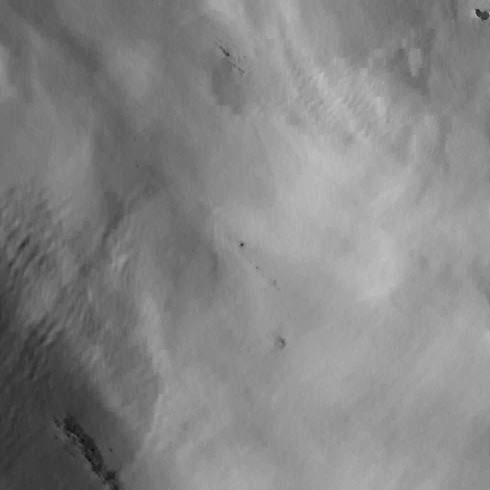}
				\includegraphics[width=.23\textwidth, angle =90]{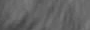}\\[1.5ex]
				 M1 & M2 \\
		\includegraphics[width=.23\textwidth, angle =90]{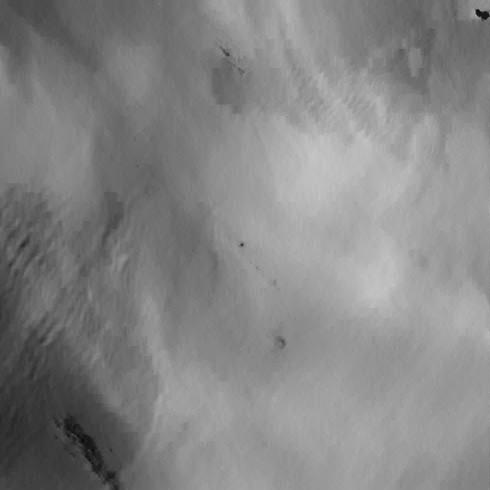}
				\includegraphics[width=.23\textwidth, angle =90]{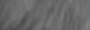}\hspace{0.1cm}
		&\includegraphics[width=.23\textwidth, angle =90]{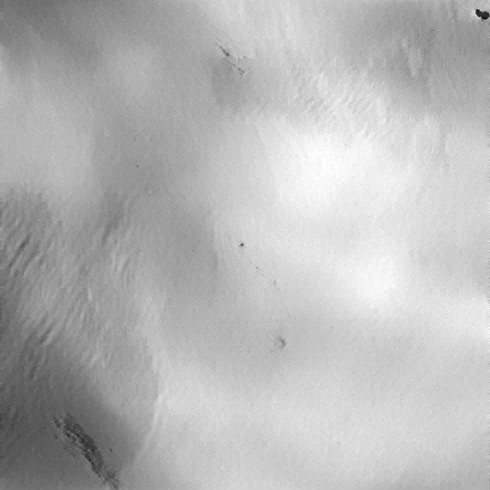}
				\includegraphics[width=.23\textwidth, angle =90]{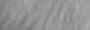}
		\end{tabular}
		\caption{Results for destriping 2D MODIS data using three different methods:
				our proposed method \protect\eqref{IC-formulation} with 
				$(\mu_1,\mu_2,\mu_3) = (0.5,1,4)$, 
				the	destriping method M1 \protect\cite{CFYL13} in \protect\eqref{model_4} with 
				$(\nu_1,\nu_2,\nu_3) = (\frac32,\frac34,1)$ and the denoising method M2 
				\protect\cite{FWL12} with $p=1, \alpha=1$.}\label{fig:2D_destr}
\end{figure*}

Our final experiments with 2D MODIS data 
(\url{http://ladsweb.nascom.nasa.gov/index.html}, also used in \cite{BL11}) 
show that our method works also fine for destriping only.
The parameters were chosen to give the best visual impression.
The results of the various methods are presented in Fig.~\ref{fig:2D_destr}.
For a comparison we have used the full model with the framelet 
regularization 
\eqref{model_4}, proposed by \cite{CFYL13}, and the method M2.
In both examples, all methods lead to visually good results.
This shows that the proposed method can also be applied for destriping, 
although it is designed for the more general setting of stripes in combination with laminar corruptions.
	
\section{Conclusions} \label{sec:conc}
We proposed a variational model 
for removing the curtaining effect in FIB tomography images.
Due to the special structure of the corruptions, 
we chose an infimal convolution model consisting of directional TV terms and also second order terms for $u$.
Alternatively, a full 3D TV term for $u$ may be applied.
The numerical experiments show that the proposed model leads to very good results. 
So far, all FIB tomography experiments share the same model parameters. 
This makes the method highly suitable for applications.
Our restoration method enables us to use the FIB tomography results for the further (statistical) analysis of 
the 3D structure of the material, which was up to now only possible for a considerable smaller area of the data.

For future work we plan to parallelize the algorithm taking also spatial domain decompositions into account.
The aim is an implementation of our approach in a way that 
it can be run simultaneously to the image acquisition process. 
Here we will also take other algorithms with potential on multicore architecture as e.g.~\cite{BCPP11,PCP11} into account.
\appendix
\section{Appendix} \label{sec:app}
{\it Proof of Proposition \ref{exist}:}
We rewrite  \eqref{IC-formulation} as
\begin{align}
		E(u,s,l) \coloneqq &\varphi_1(u) + \varphi_2(s) + 
				\varphi_3(l) + \iota_{\{(u,s,l) : u+s+l=f\}}(u,s,l). 
	\end{align}
		A feasible point is obviously given by 
		$(u,s,l)=(f,0,0)$ and $E(f,0,0) \le C$ for some
		$C \ge 0$.
		
		Assume that the infimum is not attained. Then there  exists an 
		infimal sequence $\{x_i\}_{i \in \N}$, where $x_i \coloneqq (u_i,s_i,l_i)$ so that
		\begin{align}
			 E(x_i) \le E(x_{i-1}) \le C,\quad
			 \|x_i^\tT\|_\infty \rightarrow \infty,\quad
			 E(x_i) \rightarrow \inf_{x} E(x)
		\end{align}
		as $i \rightarrow \infty$.
		We will show that we can find another sequence 
		$\{\tilde x_i\}_{i \in \N}$ with
		\begin{align}
			\|(\tilde x_i)^\tT\|_\infty 
			\le C,\quad
			E(x_i) = E(\tilde x_i).
		\end{align}
		Since this new sequence is bounded, there exists a converging subsequence.
		By the lower semi-continuity of $E$ this implies that the limit of this subsequence 
		is a minimizer of $E$ contrary to the assumption.
		
		It remains to show that there exists $\{\tilde x_i\}_{i \in \N}$.
		Consider an arbitrary element $x = (u,s,l)$ of $\{x_i\}_{i \in \N}$.
		Let the minimal and maximal value of the $j$-th $x$-$y$-plane of $l$ be 
		denoted by $m_j$, respectively $M_j$.
		Then we have by $E(x_i) \le C$ and $\mu_3 >0$ that $M_j-m_j \le C$.
		Now define 
		\begin{align}
			\tilde l = l - \begin{pmatrix} m_1 \\ \vdots \\ m_{n_z} 
			\end{pmatrix} 
			\otimes \mathbf{1}_{n_x n_y} \quad
		\mathrm{and} \quad
			\tilde s = s + \begin{pmatrix} m_1 \\ \vdots \\ m_{n_z} 
			\end{pmatrix} 
						\otimes \mathbf{1}_{n_x n_y},
		\end{align} 
		where $\mathbf{1}_{n_x n_y} = (1, \ldots, 1)^\tT \in \R^{n_x n_y}$.
		Obviously, $E(u,\tilde l, \tilde s) = E(u, l, s)$.
		Further we have $\|\tilde l\|_\infty \le C$ and by the constraint 
		on $u$ also $\|u\|_\infty \le 1$.
		Since the given image $f$ is also bounded, we conclude by
		$f=u+s+l$ that $\tilde s$ is bounded, too. This finishes the proof.	\hfill $\Box$	
\\[2ex]
{\bf Acknowledgement:}
The authors thank 
the anonymous referee for bringing model \eqref{IC-formulation_rev} to our attention,
P.~Weiss (University of Toulouse) for providing the results for the MODIS data computed by the model M2 from \cite{FWL12}, 
K.~Schladitz (Fraunhofer ITWM, Kaiserslautern), G.~Steidl, F.~Balle and T.~Beck (University of Kaiserslautern) for fruitful discussions.

We are thankful to T.~L\"ober (Nano Structuring Center Kaiserslautern) for creating FIB/SEM data of the aluminum matrix composite.
We are grateful to the German Copper Institute (Deutsches Kupferinstitut
Berufsverband e.V., DKI, D\"usseldorf, Germany) for providing bronze
samples and M.~Engstler (Material Engineering Center Saarland
(MECS), Steinbeis Research Center, Saarbr\"ucken, Germany) for providing
FIB/SEM data of bronze samples.

Funding by the German Research Foundation (DFG) 
with\-in the Research Training Group 1932 "Stochastic Models for Innovations in the Engineering Sciences", project area P3, 
is gratefully acknowledged.

\bibliography{refs_infconv}
\bibliographystyle{abbrv}

\end{document}